\documentclass[twocolumn]{autart}    
\usepackage{natbib}
\usepackage{amssymb}
\usepackage{amsmath}
\usepackage{graphicx}
\usepackage{comment}
\usepackage{booktabs}
\usepackage{bm}
\usepackage{graphicx}
\usepackage{epstopdf}
\usepackage{amssymb}
\usepackage{wasysym}
\usepackage{float}
\usepackage{color}
\usepackage{dsfont}
\usepackage{algorithm,algorithmicx}
\usepackage{algpseudocode}
\usepackage{epstopdf}
\usepackage{tabularx}
\usepackage{enumitem}
\usepackage{marginnote}
\usepackage{graphicx}          
\newcommand{\eod}{\ensuremath{\hfill\Box}}
\begin{document}

\begin{frontmatter}

\title{Event-triggered Partitioning for Non-centralized Predictive-Control-based Economic Dispatch of Interconnected Microgrids: Technical Report\thanksref{footnoteinfo} \vspace{-20pt}} 

\thanks[footnoteinfo]{
Corresponding author W.~Ananduta.}

\author[delft]{Wicak Ananduta}\ead{w.ananduta@tudelft.nl},    
\author[iri]{Carlos Ocampo-Martinez}\ead{carlos.ocampo@upc.edu}  
\address[delft]{Delft Center for Systems and Control, TU Delft}
\address[iri]{Automatic Control Department, Universitat Polit\`{e}cnica de Catalunya}  


\begin{abstract}                          
A non-centralized model predictive control (MPC) scheme for solving an economic dispatch problem of electrical networks is proposed in this paper. The scheme consists of two parts. The first part is an event-triggered repartitioning method that splits the network into a fixed number of non-overlapping sub-systems {(microgrids)}. The objective of the repartitioning procedure is to obtain self-sufficient microgrids, i.e., those that can meet their local loads using their own generation units. However, since the algorithm does not guarantee that all the resulting microgrids are self-sufficient, the microgrids that are not self-sufficient must then form a coalition with some of their neighboring microgrids. This process becomes the second part of the scheme. By performing the coalition formation, we can decompose the economic dispatch problem of the network into coalition-based sub-problems such that each subproblem is feasible. Furthermore, we also show that the solution obtained by solving the coalition-based sub-problems is a feasible but sub-optimal solution to the centralized problem. Additionally, some numerical simulations are also carried out to show the effectiveness of the proposed method.  
\end{abstract}

\end{frontmatter}

\section{Introduction}

\vspace{-10pt}
Considering the current trend and development \cite{molzahn2017,morstyn2016}, future electrical energy networks would have a high number of distributed generation and storage units. In this regard, in the energy management level, a non-centralized control scheme has been considered as a more suitable scheme than the conventional centralized one, due to the ability to deal with high computational requirement, flexibility, reliability and scalability of the non-centralized scheme \cite{molzahn2017,morstyn2016}. On the other hand, non-dispatchable generation units introduce additional uncertainty on top of the already uncertain loads. At the same time, storage units have slow dynamics that must be taken into account when solving the economic dispatch problem. In this regard, model predictive control (MPC) framework, accounting with the receding horizon principle, has been proposed to be implemented as a control scheme to the energy management of such energy systems \cite{baker2016,wang2015}.   As discussed in \cite{parisio2017,zhu2014}, the MPC framework allows us to handle components with dynamics, constraints (of physical and operational nature), and uncertainties better than traditional economic dispatch schemes.  

\vspace{-5pt}
In this paper, we discuss an energy management problem of such networks with distributed generation and storage units. In particular, we solve an MPC-based economic dispatch problem with a non-centralized scheme. We consider the microgrid framework \cite{schwaegerl2013} in which the energy network is partitioned into a group of interconnected microgrids \cite{arefifar2012,barani2018}. Each microgrid is a cluster of storage units, distributed generation units, and loads \cite{arefifar2012}. Furthermore, depending on the physical connection, it can also exchange power with the other microgrids and the main grid. More importantly, each microgrid is an independent entity that can manage itself, i.e., it has its own local controller. Therefore, in a non-centralized scheme, these microgrids cooperatively solve the economic dispatch problem of the network.

A typical non-centralized approach to solving such problems is by using a distributed optimization algorithm \cite{baker2016,wang2015,kraning2014,braun2016,guo2016} (see \cite{molzahn2017,morstyn2016} for a survey).  Such algorithms are usually iterative and require high information flow, i.e., at each iteration, each local controller must exchange information with its neighbors, with the advantage of obtaining an optimal solution. In this paper, we propose an alternative non-centralized scheme with low information flow. There are two main ingredients of the approach that we propose. The first ingredient is a proper partitioning of the network and the second ingredient is the formulation of coalition-based sub-problems, which requires a coalition formation algorithm. 

In the first part of the method, we (re)-partition the network into a fixed number of microgrids. The objective of the repartitioning scheme is to obtain self-sufficient and efficient microgrids. Roughly speaking, we consider that a microgrid is self-sufficient when it can provide its loads using its local distributed generation units.  When this goal is achieved, each microgrid does not need to rely on the other microgrids, implying a local economic dispatch problem can be solved by the controller of each microgrid. In addition, the efficiency criterion is in line with the objective of the economic dispatch problem. 
Therefore, we propose a repartitioning procedure that has low computational burden and is performed in a distributed manner. 
The repartitioning procedure that we propose is closely related to the partitioning methods presented in \cite{ananduta2019a,ananduta2019,julian2019}. The main idea of the repartitioning procedure is to move some nodes from one partition to another in order to improve the cost that has been defined. In addition, in our scheme, we consider an event-triggered mechanism, i.e., the network is only repartitioned when the event at which at least one microgrid that is not self-sufficient occurs. To the best of our knowledge, an event-triggered  repartitioning scheme has not been proposed in the literature so far.  Note that most system partitioning methods that have been proposed, e.g., \cite{fjallstrom1998,guo2016,julian2017,ocampo2011}, are intended to be implemented offline prior to the operation of the system and in a centralized fashion, whereas this paper shows how online repartitioning can be exploited to design a non-centralized control scheme and might be performed in a distributed fashion.  \color{black}

In the second part of the method, we decompose the economic dispatch problem into coalition-based sub-problems. Since the repartitioning procedure does not guarantee that the resulting microgrids are self-sufficient, each microgrid that is not self-sufficient is grouped together with some of its neighbors to form a self-sufficient coalition. 
In this regard, we propose a coalition formation procedure, which is also carried out in a distributed manner. Furthermore, coalition-based economic dispatch sub-problems are formulated. These problems are solved by the local controllers of the microgrids in order to obtain a feasible but possibly sub-optimal solution to the centralized economic dispatch problem. 
The coalition-based economic dispatch approach is inspired by the coalitional control framework \cite{fele2017,muros2017,fele2018}. In the coalitional control scheme, the sub-systems are clustered into several coalitions based on the relevance of the agents, e.g., the degree of coupling among the sub-systems. Furthermore, the sub-systems that are in the same coalition cooperatively compute their control inputs. 
 In our method, the coalitions formed are based on the necessity to maintain the feasibility of the economic dispatch sub-problems, which can be perceived as the relevance of the agents to the problem itself. \color{black} 


To summarize, the main contribution of this paper is a novel non-centralized economic dispatch method for electrical energy networks with distributed generation and storage units. The methodology has less intensive communication flows compared to typical distributed optimization methods, thus suitable with online optimization of the MPC framework. To that end, the methodology combines an event-triggered repartitioning approach with the aim of producing self-sufficient and efficient microgrids and a procedure to form self-sufficient coalitions of microgrids to solve the economic dispatch problem. This paper also provides the analysis of the proposed methodology, including the outcomes of the repartitioning and coalition formation algorithms, as well as an upper bound for the suboptimality of the proposed scheme. 
The methodology that we present in this paper is an extension of that in \cite{ananduta2019a}, where a periodical repartitioning scheme for a fully decentralized scheme is proposed. Additionally, a feasibility issue, which arises from microgrids that are not self-sufficient and can be found when using the scheme in \cite{ananduta2019a}, is solved by the coalition-based approach proposed in this paper.

The remainder of the paper is structured as follows. In Section \ref{sec:prob_form}, we provide the mathematical formulation of the economic dispatch problem and outline the proposed method. Section \ref{sec:rep_sch} presents the proposed event-triggered repartitioning scheme. In Section \ref{sec:nc_ed}, the coalition-based economic dispatch method is explained. Afterward, Section \ref{sec:opt} provides a discussion about the trade-off between the sub-optimality and communication complexity of the method. In Section \ref{sec:case_st}, a numerical study of a well-known benchmark case is presented to show the effectiveness of the proposed scheme. Finally, Section \ref{sec:concl} provides some concluding remarks and discusses future work. 

\vspace{-5pt}
\subsection*{Notations} 
\vspace{-10pt}
The sets of real numbers and integers are denoted by $\mathbb{R}$ and $\mathbb{Z}$, respectively. Moreover, for $a\in\mathbb{R}$, $\mathbb{R}_{\geq a}$ denotes all real numbers in the set \{$b: b\geq a, \ b \in \mathbb{R}$\}. A similar definition can be used for $\mathbb{Z}_{\geq a}$ and the strict inequality case. 
The set cardinality is denoted by $|\cdot|$. Finally, discrete-time instants are denoted by the subscript $k$.

\section{Problem formulation}
\label{sec:prob_form}


Let a distribution power network be represented by the undirected and connected graph $\mathcal{G}=(\mathcal{N},\mathcal{E})$, where the set of busses is denoted by $\mathcal{N}=\{1,2,\dots,n\}$ and the set of edges that connect the busses is denoted by $\mathcal{E}$, i.e., $\mathcal{E}=\{(i,j):i,j\in\mathcal{N} \}\subseteq\mathcal{N}\times\mathcal{N}$. Furthermore, the set of neighbor busses of bus $i$ is denoted by $\mathcal{N}_i$, i.e., $\mathcal{N}_i = \{j:(i,j)\in \mathcal{E} \}$. Each bus $i$ might contain an aggregate load (power demand), dispatchable or nondispatchable distributed generation units, and energy storage devices. Each of these components has operational constraints, which are assumed to be polyhedral and compact. Furthermore, each bus $i \in \mathcal{N}$ also considers power balance equations that must be satisfied  at each iteration $k$, as follows: 
\begin{align}
u^{\mathrm{g}}_{i,k} + u^{\mathrm{st}}_{i,k} + u^{\mathrm{m}}_{i,k} + \sum_{j \in \mathcal{N}_i} v^{j}_{i,k} - {d}_{i,k}=0,  \label{eq:lo_pb}\\
v^{j}_{i,k} + v^{i}_{j,k} = 0, \quad \forall j \in \mathcal{N}_i, \label{eq:coup_pb}
\end{align}
where  ${u}_{i,k}^{\mathrm{g}}\in \mathbb{R}$, $u^{\mathrm{st}}_{i,k}\in \mathbb{R}$, and $u_{i,k}^{\mathrm{m}} \in \mathbb{R}_{\geq 0}$ denote the power generated from dispatchable unit, delivered from/to storage unit, and imported from the main grid if connected, respectively.  
Furthermore, $d_{i,k}\in \mathbb{R}$ denotes the difference between uncertain loads and the power generated from non-dispatchable units, which is also uncertain. Additionally, $v_{i,k}^{j} \in \mathbb{R}$ denotes the power transferred to/from the neighbor bus $j  \in \mathcal{N}_i$. Equation \eqref{eq:lo_pb} can be considered as a local power balance, whereas \eqref{eq:coup_pb} couples two neighboring busses.

The variable $d_{i,k}$ is considered to be uncertain and its disturbance is bounded, i.e.,
\begin{align}
{d}_{i,k}&=\hat{d}_{i,k} + w_{i,k}^{\mathrm{d}}, \quad \forall i \in \mathcal{N}, \label{eq:ul}\\
|w_{i,k}^{\mathrm{d}}| &\leq \bar{w}_{i}^{\mathrm{d}}, \quad \forall i \in \mathcal{N}, \label{eq:w_ul}
\end{align}
where  $\hat{d}_{i,k}$ denotes the forecast of $d_{i,k}$ whereas  $w_{i,k}^{\mathrm{d}} \in \mathbb{R}$ and $\bar{w}_{i}^{\mathrm{d}} \in \mathbb{R}$  represent the disturbance/uncertainty and its bound, which is assumed to be known, for simplicity and since this work does not focus on handling uncertainties. A stochastic method such as the ones presented in \cite{margellos2014,ananduta2019b} can be considered as an extension. 



To state the MPC-based economic dispatch problem, define the vectors of decision variables of each bus $i\in \mathcal{N}$, which correspond to active components, by $\bm{u}_{i,k} = \begin{bmatrix}
u^{\mathrm{st}}_{i,k} & u^{\mathrm{g}}_{i,k}  & u^{\mathrm{m}}_{i,k}
\end{bmatrix}^{\top} \in \mathbb{R}^{3}$ and $\bm{v}_{i,k} = \begin{bmatrix}
v_{i,k}^{j}
\end{bmatrix}_{i\in \mathcal{N}_i}^{\top} \in \mathbb{R}^{|\mathcal{N}_i|}$. 
Furthermore, an economic quadratic cost function is considered as follows:
\begin{equation}
J_{i,k}(\bm{u}_{i,k},\bm{v}_{i,k}) = \bm{u}_{i,k}^{\top}R_i\bm{u}_{i,k}+ \bm{v}_{i,k}^{\top}Q_i\bm{v}_{i,k}, \label{eq:cost_func}
\end{equation}
where $R_i$ and $Q_i$ are positive definite diagonal matrices of suitable dimensions. Therefore, the optimization problem behind an MPC-based economic dispatch is stated as follows:
\begin{subequations}
	\begin{align}
&	\min_{\{\{(\bm{u}_{i,\ell},\bm{v}_{i,\ell})\}_{i \in \mathcal{N}}\}_{\ell=k}^{k+h-1}} \sum_{i \in \mathcal{N}} \sum_{\ell =k}^{k+h-1}  J_{i,\ell}(\bm{u}_{i,\ell},\bm{v}_{i,\ell}) \label{eq:net_cost_func}\\
& \qquad \qquad	\text{s.t. }   
	 (\bm{u}_{i,\ell},\bm{v}_{i,\ell}) \in \mathcal{P}_i, \ \forall i \in \mathcal{N}, \label{eq:net_loc_cons}\\
	&  \qquad \qquad \qquad  v_{i,\ell}^j + v_{j,\ell}^i = 0, \ \forall j\in\mathcal{N}_i, \ \ \forall i \in \mathcal{N}, \label{eq:net_coup_cons}
	\end{align}
	\label{eq:MPC_net}%
\end{subequations}
for all $\ell \in \{k,\dots,k+h-1\}$, where $h \in \mathbb{Z}_{\geq 1}$ denotes the prediction horizon. The set $\mathcal{P}_i\subset \mathbb{R}^{3+|\mathcal{N}_i|}$, for each $i\in\mathcal{N}$, is assumed to be a compact polyhedral set such that \eqref{eq:lo_pb}, \eqref{eq:ul}, and \eqref{eq:w_ul} as well as the operational constraints of the active components hold. We refer to \cite{ananduta2019a} for a more detailed description of such operational constraints. 

\color{black}
\begin{assum}\label{as:feas_ED_cent}
	For each $k \in \mathbb{Z}_{\geq 0}$, a feasible set of Problem \eqref{eq:MPC_net} exists. 
\end{assum}

\begin{rem} 	Assumption 1 is considered in order to ensure that the proposed scheme obtains a solution. In practice, the satisfaction of this assumption is achieved either if the network is connected with the main grid, which is usually assumed as an infinite source of power, or if it is not connected with the main grid, the total power that can be generated by the distributed generators is sufficiently larger than the loads within the network.
\end{rem}



\begin{figure}
	\centering
	\includegraphics[scale=0.5]{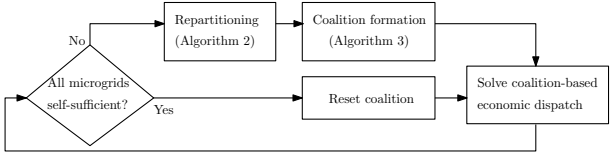}
	\caption{The overall scheme of the proposed method.
	}
	\label{fig:prop_sch}
\end{figure}

In this paper, we solve Problem \eqref{eq:MPC_net} in a non-centralized fashion, where there exists $m$ local controllers. Thus, the network must be partitioned into $m$ sub-systems, each of which is assigned to a local controller. \color{black}Then, the controllers cooperatively solve Problem \eqref{eq:MPC_net}. To that end, Problem \eqref{eq:MPC_net}, which has coupling constraints \eqref{eq:net_coup_cons}, must be decomposed. Our main idea is to decompose Problem \eqref{eq:MPC_net} into a number of sub-problems, not larger than $m$, such that each sub-problem can be solved independently. As we will show in the next sections, the independence of each sub-problem depends on the self-sufficiency of the microgrids, i.e., the ability to meet local load using local production. Therefore, we propose an event-triggered repartitioning and coalition formation procedures to obtain self-sufficient partitions. 
A flow diagram that summarizes the overall method is shown in Figure \ref{fig:prop_sch}. 


\section{Event-triggered repartitioning scheme}
\label{sec:rep_sch}
Since the loads in the network vary over time, the network might need to be repartitioned to maintain the level of self-sufficiency. In this section, first we state the repartitioning problem. Afterward, we present the repartitioning process, which includes when and how the repartitioning must be performed. Furthermore, we also provide an analysis of the  repartitioning outcome.

\subsection{Repartitioning Problem Formulation}
Prior to presenting the methodology that we propose, we establish some definitions that will be used throughout the remainder of the paper. 

\begin{defn}[Non-overlapping partition]
	\label{def:nop}
	The set $\boldsymbol{\mathcal{M}}~=~\{\mathcal{M}_1, \mathcal{M}_2, \dots, \mathcal{M}_m \}$ defines $m$ non-overlapping partitions of graph $\mathcal{G}=(\mathcal{N},\mathcal{E})$ if $\bigcup_{p=1}^m \mathcal{M}_p = \mathcal{N}$ and $\mathcal{M}_p \cap \mathcal{M}_q = \emptyset$, for any $\mathcal{M}_p,\mathcal{M}_q \in \boldsymbol{\mathcal{M}}$ and $p \neq q$. 
\end{defn}

\begin{defn}[Local imbalance]
	The local power imbalance of a subset of nodes $\mathcal{M}\subseteq\mathcal{N}$ at any $k\geq0$, denoted by $\Delta_{\mathcal{M},k}^{\mathrm{im}}$, is defined as
	\begin{equation}
		\Delta_{\mathcal{M},k}^{\mathrm{im}}= \sum_{i \in \mathcal{M}} \left(-\bar{u}_{i}^{\mathrm{g}} + {d}_{i,k}\right),
		\label{eq:p_im}
	\end{equation}
	where $\bar{u}_{i}^{\mathrm{g}}$ denotes the maximum capacity of dispatchable generation units, whereas ${d}_{i,k}$ follow \eqref{eq:ul} and the worst case disturbance, which is shown in \eqref{eq:w_ul}, is considered. \eod
	\label{def:local_imb}
\end{defn} 

\begin{defn}[Self-sufficiency]
	A subset of nodes $\mathcal{M}\subseteq\mathcal{N}$ at any $k\geq0$ is self-sufficient if it has non-positive local imbalance at any step along the prediction horizon $h$, i.e., $\Delta_{\mathcal{M},\ell}^{\mathrm{im}} \leq 0,$ for all $\ell \in \{k,k+1,\dots,k+h-1 \}$. \eod
	\label{def:ss_mg}
\end{defn}
\begin{defn}[Imbalance cost]
	The imbalance cost of microgrid $\mathcal{M}_{p,k}\in\bm{\mathcal{M}}_k$ at any $k\geq0$, denoted by $J_{p,k}^{\mathrm{im}}$, is defined as
$	J_{p,k}^{\mathrm{im}} = \sum_{\ell =k}^{k+h-1}\max \left(0,\Delta_{\mathcal{M}_{p,k},\ell}^{\mathrm{im}}\right),$
	where $\Delta_{\mathcal{M}_{p,k},\ell}^{\mathrm{im}}$ is defined based on \eqref{eq:p_im}. \eod
	\label{def:im_cost}
\end{defn}
\begin{defn}[Efficiency cost]
	\label{def:eff_cost}
	The efficiency cost of microgrid $\mathcal{M}_{p,k}\in\bm{\mathcal{M}}_k$ at any $k\geq0$, is defined  as follows:
	\begin{equation*}
	\begin{aligned}
	J_{p,k}^{\mathrm{ef}}=&\min_{\{\{\bm{u}_{i,\ell}\}_{i \in \mathcal{M}_{p,k}}\}_{\ell=k}^{k+h-1}} \sum_{\ell =k}^{k+h-1} \left(J_{p,\ell}^{\mu} + J_{p,\ell}^{{\epsilon}}\right)\\
	& \ \quad \mathrm{s.t.} \
	(\bm{u}_{i,\ell},\bm{v}_{i,\ell}) \in \mathcal{P}_i, \ \forall i \in \mathcal{M}_{p,k}, \\
	& \quad \quad v_{i,\ell}^j + v_{j,\ell}^i = 0, \ \forall j\in\mathcal{N}_i\cap\mathcal{M}_{p,k}, \ \forall i \in \mathcal{M}_{p,k},\\ 
	& \quad \quad \forall\ell \in\{k,\dots, k+h-1 \},
	\end{aligned}
	\end{equation*}
	where $J_{p,\ell}^{{\epsilon}}$ adds extra cost on the power transferred between one microgrid to another in order to minimize the dependency on the neighbors and is defined by
	$J_{p,\ell}^{{\epsilon}} = \sum_{i\in \mathcal{M}_{p,k}^{\mathrm{c}}} \sum_{j \in \mathcal{N}_i\backslash\mathcal{M}_{p,k}} c_{i}^{\mathrm{et}}(p_{ij,\ell}^{\mathrm{t}})^2,$
	where $c_{i}^{\mathrm{et}} \in \mathbb{R}_{\geq 0}$ is the extra per-unit cost of transferring power.\eod
\end{defn}

Now, we state the repartitioning problem that will be solved supposing that the network is triggered to perform the repartitioning. 
First, assume that the network is initially partitioned into $m$ non-overlapping microgrids and denote the set of initial partition at $k=0$ by $\boldsymbol{\mathcal{M}}_0~=~\{\mathcal{M}_{1,0}, \mathcal{M}_{2,0}, \dots, \mathcal{M}_{m,0} \}$. 
Thus, for some time instants, at which the system must perform repartitioning,  the optimization problem that   must be solved is stated as follows: 
\begin{equation}
	\begin{aligned}
	&\min_{\boldsymbol{\mathcal{M}}_k}  \ \sum_{p=1}^m J^{\pi}(\mathcal{M}_{p,k})\\
	& \text{s.t. } \boldsymbol{\mathcal{M}}^{(0)}_k = \boldsymbol{\mathcal{M}}_{k-1},\\
	& \mathcal{M}_{p,k} \in \boldsymbol{\mathcal{M}}_k \text{ is non-overlapping and connected}.%
	\end{aligned}
	\label{eq:part_prob}%
\end{equation}
The cost function $J^{\pi}(\mathcal{M}_{p,k})$ is defined by 
\begin{equation}
J^{\pi}(\mathcal{M}_{p,k})= \alpha J_{p,k}^{\mathrm{im}} + J_{p,k}^{\mathrm{ef}}, 
\label{eq:j_pi}
\end{equation}
where $\alpha$ is the tuning parameter to determine the trade-off between both imbalance and efficiency costs. 
Moreover, each microgrid must be connected, i.e., the subgraph formed by each microgrid is connected. This constraint is imposed to avoid decoupling among the nodes within each microgrid. 
Furthermore, $\boldsymbol{\mathcal{M}}^{(0)}_k$ denotes the initial partition at time step $k$, which is equal to the partition at the previous time step, $k-1$. 
In addition, Assumption \ref{as:in_par}, which is related to the initial partition $\boldsymbol{\mathcal{M}}_0$, is considered.
\begin{assum}
	\label{as:in_par}
	The initial partition $\boldsymbol{\mathcal{M}}_0$ is non-overlapping with connected microgrids. 
\end{assum}

\subsection{Repartitioning Process}

The repartitioning process consists of two main steps. The first step is to determine whether the system must perform the repartitioning and the second step is to actually perform the repartitioning. The event that triggers a repartitioning process is the existence of a microgrid that is not self-sufficient (c.f. Definition \ref{def:ss_mg}). In this regard, the triggering mechanism is provided in Algorithm \ref{alg:trig}. 
\begin{alg}[Triggering mechanism]
	\label{alg:trig}
	  \hfill
	\begin{enumerate}
		\item For each microgrid $\mathcal{M}_{p,k-1} \in\boldsymbol{\mathcal{M}}_{k-1}$, evaluate its self-sufficiency at $k$, based on Definition \ref{def:ss_mg}.
		\item If a microgrid is not self-sufficient, raise a flag to start repartitioning procedure. Otherwise wait until all microgrids perform step 1.
		\item If the flag to start the repartitioning procedure is not raised, then keep the  current partition, i.e., $\boldsymbol{\mathcal{M}}_{k}=\boldsymbol{\mathcal{M}}_{k-1}$. \eod
	\end{enumerate}
\end{alg}

Now, we discuss the repartitioning process, where the controllers cooperatively solve Problem \eqref{eq:part_prob}. We propose an iterative local improvement algorithm that can be performed in a distributed and synchronous manner. 
Consider the initial partition $\boldsymbol{\mathcal{M}}_k^{(0)}$. Moreover, denote the iteration number by superscript $(r)$ and the set of boundary busses of microgrid $\mathcal{M}_{p,k}$, i.e., busses that are connected (coupled) to at least one bus that belongs to another microgrid by $\mathcal{M}_{p,k}^{\mathrm{c}}=\{i: (i,j) \in \mathcal{E}, i\in \mathcal{M}_{p,k}, j\in \mathcal{N}\backslash\mathcal{M}_{p,k} \}$. Then, the repartitioning procedure is stated in Algorithm \ref{alg:repartitioning}. 
Note that Algorithm \ref{alg:repartitioning} can be stopped when it reaches a predefined maximum  number of iteration $\bar{r}$.  \color{black}

\begin{alg}[Repartitioning procedure]
	\label{alg:repartitioning}
	Suppose that microgrid $\mathcal{M}_{p,k}$ is chosen to propose a node that will be moved at the $r^{\text{th}}$ iteration. Then, the steps at each iteration are described below:
	
	\begin{enumerate}
		\item Microgrid $\mathcal{M}_{p,k}$ computes $J^{\pi}(\mathcal{M}_{p,k}^{(r)})$, which is the local cost function at the $r^{\text{th}}$ iteration, based on \eqref{eq:j_pi}.
		\item Microgrid $\mathcal{M}_{p,k}$ computes a node that will be offered to be moved, denoted by $\theta_p$, as follows:
		\begin{equation}
		\theta_p \in \arg\min_{\theta \in \mathcal{M}_{p,k}^{\theta(r)}} J^{\pi}(\mathcal{M}_{p,k}^{(r)}\backslash\{\theta \}),
		\label{eq:theta_p}
		\end{equation}
		where $\mathcal{M}_{p,k}^{\theta(r)} \subseteq \mathcal{M}_{p,k}^{\mathrm{c}(r)}$ is a subset of boundary busses that do not disconnect microgrid $\mathcal{M}_{p,k}$ when removed, i.e., the graph form by $\mathcal{M}_{p,k}^{(r)} \backslash \{\theta \}$, for $\theta \in \mathcal{M}_{p,k}^{\theta(r)}$, is connected. The node   $\theta_p$ is randomly selected from the set of minimizers of \eqref{eq:theta_p}.
		\item 
		Microgrid $\mathcal{M}_{p,k}$ computes the local cost difference if $\theta_p$ is moved out from microgrid $\mathcal{M}_{p,k}$, i.e.,
		\begin{equation}
		\Delta J^{\pi(r)}_p = J^{\pi}(\mathcal{M}_{p,k}^{(r)}\backslash\{\theta_p \}) - J^{\pi}(\mathcal{M}_{p,k}^{(r)}).
		\label{eq:delta_Jp}
		\end{equation} 
		\item Microgrid $\mathcal{M}_{p,k}$ shares the information of $\theta_p$ and $\Delta J^{\pi(r)}_p$ to the related neighboring microgrids $\mathcal{M}_{q,k}^{(r)} \in \mathcal{N}_{\theta_p}'=
		\{\mathcal{M}_{q,k}^{(r)}: (\theta_p,j) \in \mathcal{E}, j \in \mathcal{M}_{q,k}^{(r)} \}$.
		\item  All neighbors $\mathcal{M}_{q,k}^{(r)} \in \mathcal{N}_{\theta_p}'$ compute the expected total cost difference if $\theta_p$ is moved from microgrid $\mathcal{M}_{p,k}^{(r)}$ to microgrid $\mathcal{M}_{q,k}^{(r)}$, as follows:
		\begin{equation}
		\Delta J^{\pi(r)}_{\mathrm{t},q} = J^{\pi}(\mathcal{M}_{q,k}^{(r)}\cup\{\theta_p \}) - J^{\pi}(\mathcal{M}_{q,k}^{(r)})+\Delta J^{\pi(r)}_p,
		\label{eq:delta_Jq}
		\end{equation}
		and send the information of $\Delta J^{\pi(r)}_{\mathrm{t},q}$ to microgrid $\mathcal{M}_{p,k}$.
		\item Microgrid $\mathcal{M}_{p,k}$ selects the neighbor that will receive $\theta_p$ as follows:
		$q^{\star} \in \arg\min_{q \in \mathcal{N}_{\theta_p}'}\Delta J^{\pi(r)}_{\mathrm{t},q},
		$
		where $q^{\star}$ is randomly chosen from the set of minimizers.
		\item  If $\Delta J^{\pi(r)}_{\mathrm{t},q^{\star}} \leq 0$, then the partition is updated as follows:
		$\mathcal{M}_{p,k}^{(r+1)} = \mathcal{M}_{p,k}^{(r)}\backslash \{\theta_p\},\
		\mathcal{M}_{q^{\star},k}^{(r+1)} = \mathcal{M}_{q^{\star},k}^{(r)}\cup \{\theta_p\}.
		$ 
		Otherwise, the algorithm jumps to the next iteration, $r+1$. \eod
	\end{enumerate}
\end{alg}

\begin{prop}
	\label{prop:sol_alg}
	Let $\boldsymbol{\mathcal{M}}_0$ be the initial partition at $k=0$ and Assumption \ref{as:in_par} hold. At any time instant at which the repartitioning process is triggered, the output of Algorithm \ref{alg:repartitioning} is a non-overlapping partition with connected microgrids and converges to a local minimum. \eod 
\end{prop}

\begin{pf*}{Proof.}
	Define by $\kappa$ the time instant at which the repartitioning process is triggered, i.e., there exists at least one microgrid in $\boldsymbol{\mathcal{M}}_{\kappa}$ that is not self-sufficient. Let $\kappa_0$ be the first (smallest) repartitioning instant. Notice that the initial partition $\boldsymbol{\mathcal{M}}^{(0)}_{\kappa}$, at any repartitioning instant $\kappa > \kappa_0$, equals to the solution of Algorithm \ref{alg:repartitioning} at the previous repartitioning instant. Therefore, if at $\kappa_0$ the assertion holds, then it also holds for any repartitioning instants. Hence, now we only need to evaluate the outcome of the repartitioning process at $\kappa_0$. 
	Since the system is not repartitioned when $k<\kappa_0$, the initial partition at $\kappa_0$, $\boldsymbol{\mathcal{M}}^{(0)}_{\kappa_0}=\boldsymbol{\mathcal{M}}_0$, is non-overlapping with connected microgrids due to Assumption \ref{as:in_par}. Moreover, at any iteration of the repartitioning procedure, the node proposed to be moved is selected from $\mathcal{M}_{p,{\kappa_0}}^{\theta (r)}$, which is the set of boundary nodes that do not cause the disconnection of the associated microgrid when removed (see \eqref{eq:theta_p}). At the end of the iteration, either one node is moved from one microgrid to another or no node is moved. These facts imply that, at the end of any iteration, $\boldsymbol{\mathcal{M}}^{(r)}_{\kappa_0}$ remains non-overlapping and the connectivity of each microgrid is maintained. 
	Now, we show the convergence of the repartitioning solution. To this end, let the total cost at the beginning of iteration $r$  be denoted by $J^{\pi(r)}_{\mathrm{t}} = \sum_{p=1}^m J^{\pi}(\mathcal{M}_{p,{\kappa}})$. The convergence is proved by showing that the evolution of the total cost is non-increasing. Suppose that $\theta_p$ is moved from microgrid $p$ to microgrid $q^{\star}$. Therefore, we have
	\vspace{-5pt}
	\begin{align*}
	J^{\pi(r+1)}_{\mathrm{t}}- J^{\pi(r)}_{\mathrm{t}} 
	&= J^{\pi}(\mathcal{M}_{p,k}^{(r+1)}) - J^{\pi}(\mathcal{M}_{p,{\kappa}}^{(r)}) \\& \quad + J^{\pi}(\mathcal{M}_{q^{\star},{\kappa}}^{(r+1)}) - J^{\pi}(\mathcal{M}_{q^{\star},{\kappa}}^{(r)})\\
	&=\Delta J^{\pi(r)}_{\mathrm{t},q^{\star}} \leq 0.
	\end{align*}
	The first equality follows from the fact that only the local costs of microgrids $p$ and $q^{\star}$ change after iteration $r$. The second equality follows from \eqref{eq:delta_Jp} and \eqref{eq:delta_Jq}, and the inequality comes from the condition imposed in step 7 of Algorithm \ref{alg:repartitioning}, where $\theta_p$ is not moved if  $\Delta J^{\pi(r)}_{\mathrm{t},q^{\star}} > 0$. When no node is moved, $J^{\pi(r+1)}_{\mathrm{t}}- J^{\pi(r)}_{\mathrm{t}}=0$. \eod
\end{pf*}


\section{Coalition-based economic dispatch scheme}
\label{sec:nc_ed}
In this section, the non-centralized economic dispatch scheme based on the previously explained repartitioning approach is discussed. 
We let each self-sufficient microgrid to solve its local economic dispatch problem and does not allow this microgrid to exchange power with its neighbors. Therefore, by imposing an additional constraint, self-sufficient microgrids do not need to communicate with its neighbors to dispatch its components. However, a fully decentralized method can only be performed if all microgrids are self-sufficient. For any microgrid that is not self-sufficient, its local economic dispatch problem might be infeasible since local power production is not enough to meet the load.  Since the repartitioning outcome does not guarantee the self-sufficiency of each microgrid, then the microgrids that are not self-sufficient form a coalition with some other microgrids such that the resulting economic dispatch problem that must be solved by each coalition is feasible. Note that in general, Problem \eqref{eq:MPC_net} might actually have feasible solutions that require high power exchange, implying it might be impossible to partition the network into self-sufficient microgrids.

\subsection{Coalition formation}
In order to describe the coalition formation procedure, denote by $\mathcal{C}_{p,k}$ and $\mathcal{D}_{p,k}$ the set of nodes and the set of microgrids that belong to coalition $p$, respectively. We  assign one pair $(\mathcal{C}_{p,k},\mathcal{D}_{p,k})$ to each microgrid $\mathcal{M}_{p,k}$ to keep tracking the nodes and neighboring microgrids with which it forms a coalition. The coalition formation mechanism is described in Algorithm \ref{alg:clust}. 
\begin{alg}[Coalition formation procedure]
	\label{alg:clust}
	\hfill
	
	Each microgrid $\mathcal{M}_{p,k}$ defines $\mathcal{C}_{p,k}^{(0)}=\mathcal{M}_{p,k}$ and $\mathcal{D}_{p,k}^{(0)}=\{\mathcal{M}_{p,k}\}$. 
	While $r < m-1$, do:
	\begin{enumerate}
		
		\item Each microgrid $\mathcal{M}_{p,k}$ evaluates whether its coalition is self-sufficient based on Definition \ref{def:ss_mg}, i.e., whether $ \Delta_{\mathcal{C}_{p,k}^{(r)},\ell}^{\mathrm{im}} \leq 0, \quad \forall \ell \in \{k,\dots,k+h-1 \},$
		holds true. 
		\item If coalition $\mathcal{C}_{p,k}^{(r)}$ is self-sufficient, then microgrid $\mathcal{M}_{p,k}$ waits for commands until $r=m-1$.  
		\item Otherwise, microgrid $\mathcal{M}_{p,k}$ initiates a coalition merger by sending $\Delta_{\mathcal{C}_{p,k}^{(r)},\ell}^{\mathrm{im}}$, for all $\ell \in \{k,\dots,k+h-1 \}$, to the microgrids that do not belong to coalition $\mathcal{C}_{p,k}^{(r)}$ but they have physical connections with at least one bus in coalition $\mathcal{C}_{p,k}^{(r)}$, i.e.,  $\mathcal{M}_{q,k} \in \mathcal{N}_{p,k}^{\mathrm{c}}= \{\mathcal{M}_{q,k}: (i,j)\in \mathcal{E}, i \in \mathcal{C}_{p,k}^{(r)}, j \in \mathcal{M}_{q,k}, \mathcal{C}_{q,k}^{(r)} \neq \mathcal{C}_{p,k}^{(r)}  \}$. Note that if $\mathcal{N}_{p,k}^{\mathrm{c}}=\emptyset$, then $\mathcal{M}_{p,k}$ waits for commands until  $r=m-1$. \color{black}
		\item For each neighbor $\mathcal{M}_{q,k} \in \mathcal{N}_{p,k}^{\mathrm{c}}$, if it is not communicating with another microgrid, then it computes 
		$J_{q}^{\mathrm{cim}}= \sum_{\ell =k}^{k+h-1} \max\left(0,\Delta_{\mathcal{C}_{q,k}^{(r)},\ell}^{\mathrm{im}}+\Delta_{\mathcal{C}_{p,k}^{(r)},\ell}^{\mathrm{im}}\right).$ Otherwise, $J_{q}^{\mathrm{cim}}=\infty$. Then, it sends back $J_{q}^{\mathrm{cim}}$ to coalition $\mathcal{C}_{p,k}^{(r)}$.  
		\item Based on $J_{q}^{\mathrm{cim}}$, microgrid $\mathcal{M}_{p,k}$ chooses the neighbor that it will merge with as a coalition, as follows: 
		\begin{equation*}
		q^{\star}\in\arg\min_{\mathcal{M}_{q,k} \in \mathcal{N}_{p,k}^{\mathrm{c}}} J_{q}^{\mathrm{cim}} \quad \mathrm{s.t.} \ J_{q}^{\mathrm{cim}} < \infty.
		\end{equation*}
		\item Update the coalition sets, i.e.,  $\mathcal{C}_{\rho,k}^{(r+1)} = \mathcal{C}_{\rho,k}^{(r)}\cup \mathcal{C}_{q^{\star},k}^{(r)}$ and $\mathcal{D}_{\rho,k}^{(r+1)} = \mathcal{D}_{\rho,k}^{(r)} \cup \mathcal{D}_{q^{\star},k}^{(r)}$, for all $\mathcal{M}_{\rho,k}\in\mathcal{D}_{p,k}^{(r)}$ and $\mathcal{C}_{\rho,k}^{(r+1)} = \mathcal{C}_{\rho,k}^{(r)}\cup \mathcal{C}_{p,k}^{(r)}$ and $\mathcal{D}_{\rho,k}^{(r+1)} = \mathcal{D}_{\rho,k}^{(r)} \cup \mathcal{D}_{p,k}^{(r)}$, for all $\mathcal{M}_{\rho,k}\in\mathcal{D}_{q^{\star},k}^{(r)}$.
		\item $r\leftarrow r+1$ and go back to step 1. 
	\end{enumerate}
\end{alg}
\begin{prop}
	\label{prop:alg_clus}
	By performing Algorithm \ref{alg:clust}, either all resulting coalitions $\mathcal{C}_{p,k}^{(m-1)}$, for $p=1,\dots,m$, are self-sufficient or all coalitions are merged, i.e., $\mathcal{C}_{p,k}^{(m-1)}=\mathcal{N}$, for $p=1,\dots,m$. \eod
\end{prop}


\begin{pf*}{Proof.}
	 At each iteration $r<m-1$, the evaluation in step 1 has two mutually exclusive outcomes:
	 \begin{enumerate}
	 		 	\setlength\itemsep{0.0em}
	 	\item All coalitions are self-sufficient.
	 	\item There exist some coalitions that are not self-sufficient.
	 \end{enumerate}
	 In case 1, we have that $\mathcal{C}_{p,k}^{(m-1)}=\mathcal{C}_{p,k}^{(r)}$, for all $p=1,\dots,m$ since the coalitions do not change from the $r^{\mathrm{th}}$ iteration until the $(m-1)^{\mathrm{th}}$ iteration.  Note that when all microgrids  $\mathcal{M}_{p,k} \in\boldsymbol{\mathcal{M}}_k$ are self-sufficient, then $\mathcal{C}_{p,k}^{(0)}$, for all $p=1,\dots,m$, are self-sufficient. Therefore, this case is also included here. 
	 In case 2, according to steps 3-6, at least one of the coalitions that are not self-sufficient will be merged with one of its neighboring coalitions at the next iteration $r+1$. Since the number of initial coalitions is finite ($m$), then if case 2 keeps occurring, all coalitions will be merged, i.e., $\mathcal{C}_{p,k}=\mathcal{N}$, for all $p=1,\dots,m$, at a finite number of iterations. Otherwise, case 1 will occur. Furthermore, in case 2, the minimum number of coalitions that can perform steps 3-6 (merging with one of its neighboring coalitions) is one. If, for $r\geq1$, only one coalition merges with one of its neighbors, then it requires $m-1$ iterations to merge all coalitions. 
	  \eod \color{black}
\end{pf*}

\subsection{Non-centralized economic dispatch}
In this section, we outline the proposed scheme to solve Problem \eqref{eq:MPC_net} based on the coalitions that have been formed. Note that when all microgrids $\mathcal{M}_{p,k}, \ p=1,\dots,m,$ are self-sufficient, the coalitions are reset as in the initialization of Algorithm \ref{alg:clust}, i.e., $\mathcal{C}_{p,k}=\mathcal{M}_{p,k}$, for $p=1,\dots,m$. First, we reformulate Problem \eqref{eq:MPC_net} as shown in Proposition \ref{lem:net2co}.
\begin{prop}
	\label{lem:net2co}
	Suppose that, at time instant $k$, the network is partitioned into $m$ non-overlapping microgrids, defined by the set $\boldsymbol{\mathcal{M}}_k~=~\{\mathcal{M}_{1,k}, \mathcal{M}_{2,k}, \dots, \mathcal{M}_{m,k} \}$. Furthermore, coalitions of microgrids, denoted by $\mathcal{C}_{p,k}$, for $p=1,\dots,m$, are formed based on Algorithm \ref{alg:clust}. 
	Then, Problem \eqref{eq:MPC_net} is equivalent to
	\begin{subequations}
		\begin{align}
		&\min_{\{\{(\bm{u}_{i,\ell},\bm{v}_{i,\ell})\}_{i \in \mathcal{N}}\}_{\ell=k}^{k+h-1}}  \sum_{p=1}^{m} \sum_{i \in \mathcal{M}_{p,k}} \sum_{\ell =k}^{k+h-1}  J_{i,\ell}(\bm{u}_{i,\ell},\bm{v}_{i,\ell}) \label{eq:mg_cost_func}\\
		&\mathrm{s.t. }  
		\quad (\bm{u}_{i,\ell},\bm{v}_{i,\ell}) \in \mathcal{P}_i, \ \forall i \in \mathcal{C}_{p,k},  \label{eq:mg_loc_cons}\\
		&\qquad \quad v_{i,\ell}^j + v_{j,\ell}^i = 0, \ \forall j \in \mathcal{N}_i\cap\mathcal{C}_{p,k}, \  \forall i \in \mathcal{C}_{p,k},  \label{eq:mg_loc_cons2}\\
		&\qquad \quad v_{i,\ell}^j + v_{j,\ell}^i = 0, \ \forall j \in \mathcal{N}_i\backslash\mathcal{C}_{p,k}, \ \forall i \in \mathcal{C}_{p,k},  \label{eq:mg_coup_cons}
		\end{align}
		\label{eq:MPC_mg}%
	\end{subequations}
	for all $\mathcal{C}_{p,k}$, where $p=1,\dots,m$, and all $\ell \in\{k,\dots, k+h-1 \}$. \eod 
\end{prop}

\begin{pf*}{Proof.}
	Since $\mathcal{M}_{p,k}$, for each $p=1,\dots,m$, is non-overlapping, the cost function \eqref{eq:net_cost_func} is equal to  \eqref{eq:mg_cost_func}. Moreover, since $\bigcup_{p=1}^m \mathcal{C}_{p,k} = \mathcal{N}$,  \eqref{eq:net_loc_cons} is equivalent to \eqref{eq:mg_loc_cons}. Furthermore, \eqref{eq:net_coup_cons} is decomposed into \eqref{eq:mg_loc_cons2} and \eqref{eq:mg_coup_cons}.  \eod
\end{pf*}
\begin{rem}
	\label{rem:const_coal}
	For each coalition $\mathcal{C}_{p,k}$, \eqref{eq:mg_loc_cons} and \eqref{eq:mg_loc_cons2} are local constraints. Particularly for the constraints in \eqref{eq:mg_loc_cons2}, some of them might involve two different microgrids. Meanwhile, \eqref{eq:mg_coup_cons} are coupling constraints with other coalitions. \eod
\end{rem}
By decomposing Problem \eqref{eq:MPC_mg}, we formulate the decentralized MPC-based economic dispatch problem that must be solved at each coalition $\mathcal{C}_{p,k}$, for $p=1,\dots,m$, as follows:
\begin{subequations}
	\begin{align}
	&\min_{\{\{(\bm{u}_{i,\ell},\bm{v}_{i,\ell})\}_{i \in \mathcal{C}_{p,k}}\}_{\ell=k}^{k+h-1}} \sum_{i \in \mathcal{C}_{p,k}} \sum_{\ell =k}^{k+h-1}  J_{i,\ell}(\bm{u}_{i,\ell},\bm{v}_{i,\ell})\\
	&\text{s.t. } \quad  (\bm{u}_{i,\ell},\bm{v}_{i,\ell}) \in \mathcal{P}_i,   \label{eq:co_loc_cons}\\
		&\qquad \quad v_{i,\ell}^j + v_{j,\ell}^i = 0, \ \forall j \in \mathcal{N}_i\cap\mathcal{C}_{p,k},   \label{eq:co_loc_cons2}\\
		&\qquad \quad v_{i,\ell}^j = 0, \ \forall j \in \mathcal{N}_i\backslash\mathcal{C}_{p,k},  \label{eq:coup_dec_const}
	\end{align}
	\label{eq:lo_mpc}%
\end{subequations}
for all $i \in \mathcal{C}_{p,k}$ and $\ell \in\{k,\dots, k+h-1 \}$. Note that if microgrids $\mathcal{M}_{p,k}$ and $\mathcal{M}_{q,k}$ belong to the same coalition, i.e., $\mathcal{C}_{p,k}=\mathcal{C}_{q,k}$,  then they must cooperatively solve the same problem in a distributed manner. Additionally, if $c=m$, then each microgrid is self-sufficient, implying a fully decentralized scheme (without communication) is applied to the network. On the other hand, if $c=1$, then a fully distributed scheme (with neighbor-to-neighbor communication) is applied to the network.

Now, we show that Problem \eqref{eq:lo_mpc}, for any coalition, has a solution. Furthermore, the solution to Problem \eqref{eq:lo_mpc} is also a feasible solution to the original problem \eqref{eq:MPC_net}.
\begin{prop}
	Suppose that Assumption \ref{as:feas_ED_cent} holds and let the coalitions $\mathcal{C}_{p,k}$, for $p=1,\dots,m$, are formed by using Algorithm \ref{alg:clust}. Then,  there exists a unique solution to Problem \eqref{eq:lo_mpc}, for each coalition $\mathcal{C}_{p,k}$, where \ $p \in \{1,\dots,m\}$.
\end{prop}

\begin{pf*}{Proof.}
Since the cost function is strictly convex, the uniqueness of the solution is guaranteed provided that the feasible set is nonempty. Therefore, we only need to show that Problem \eqref{eq:lo_mpc}, for any $\mathcal{C}_{p,k}$, has a  non-empty feasible set. Consider the outcome of Algorithm \ref{alg:clust} (c.f. Proposition \ref{prop:alg_clus}). If Algorithm \ref{alg:clust} results in one coalition over the whole network, i.e., $\mathcal{C}_{p,k}=\mathcal{N}$, for $p=1,\dots,m$, then it implies that all microgrids must solve the centralized economic dispatch problem \eqref{eq:MPC_mg} cooperatively. Therefore, in this case, for any $\mathcal{C}_{p,k}$, Problem \eqref{eq:lo_mpc} is equal to Problem \eqref{eq:MPC_mg}. Due to Assumption \ref{as:feas_ED_cent}, feasible solutions to Problem \eqref{eq:MPC_mg} exist. Otherwise, Algorithm \ref{alg:clust} results in at least two different self-sufficient coalitions. In this case, we have non-positive local imbalance (c.f. Definition \ref{def:local_imb}), i.e., the worst-case uncertain imbalance between loads and non-dispatchable generation can be met cooperatively by the distributed generation units within the coalition.  Therefore, there exists a non-empty subset of feasible solution of Problem \eqref{eq:MPC_mg} such that \eqref{eq:coup_dec_const}, for all $\mathcal{C}_{p,k}$, where $p=1,\dots,m$, hold, implying the existence of a non-empty feasible set of Problem \eqref{eq:lo_mpc}. \eod
\end{pf*}
\begin{prop}
	\label{prop:feas_of_orig}
	Let $(\bm{u}_{i,\ell}^{\star},\bm{v}_{i,\ell}^{\star})$, for all $\ell \in\{k,\dots, k+h-1 \}$ and $i\in\mathcal{C}_{p,k}$, be the solution to Problem \eqref{eq:lo_mpc}, for all coalitions  $\mathcal{C}_{p,k}$, where $p=1,\dots,m$. Then, they are also a feasible solution to Problem \eqref{eq:MPC_net}. 
\end{prop}

\begin{pf*}{Proof.}
	In Proposition \ref{lem:net2co}, we show that Problem \eqref{eq:MPC_mg} is equivalent with Problem \eqref{eq:MPC_net}, therefore we only need to show that $(\bm{u}_{i,\ell}^{\star},\bm{v}_{i,\ell}^{\star})$, for all $\ell \in\{k,\dots, k+h-1 \}$, $i\in\mathcal{C}_{p,k}$, and $p=1,\dots,m$, is a feasible solution to Problem \eqref{eq:MPC_mg}. Note that Problem \eqref{eq:lo_mpc} is obtained by decomposing Problem \eqref{eq:MPC_mg}.  As can be seen, the constraints \eqref{eq:mg_loc_cons}-\eqref{eq:mg_loc_cons2} are decomposed for each coalition and considered as \eqref{eq:co_loc_cons}-\eqref{eq:co_loc_cons2} in Problem \eqref{eq:lo_mpc}. Since $(\bm{u}_{i,\ell}^{\star},\bm{v}_{i,\ell}^{\star})$, for all $\ell \in\{k,\dots, k+h-1 \}$ and $i\in\mathcal{C}_{p,k}$, satisfy the constraints \eqref{eq:co_loc_cons}-\eqref{eq:co_loc_cons2}, they also satisfy \eqref{eq:mg_loc_cons}-\eqref{eq:mg_loc_cons2}. Finally, for any $\mathcal{C}_{p,k}$, by \eqref{eq:coup_dec_const}, we know that $v_{i,\ell}^{j\star}=v_{j,\ell}^{i\star}=0$, for all $j\in\mathcal{N}_i\backslash\mathcal{C}_{p,k}$, $i\in\mathcal{C}_{p,k}$, and $\ell\in\{k,\dots, k+h-1 \}$. From this fact, we obtain that $v_{i,\ell}^{j\star}+v_{j,\ell}^{i\star}=0$ for all $j\in\mathcal{N}_i\backslash\mathcal{C}_{p,k}$ $i\in\mathcal{C}_{p,k}$, and $\ell\in\{k,\dots, k+h-1 \}$, implying the satisfaction of the constraints in \eqref{eq:mg_coup_cons}. \eod
\end{pf*}
Finally, we note that due to the following coupling constraints,
\begin{equation}
v_{i,\ell}^j + v_{j,\ell}^i = 0, \quad \forall j \in \mathcal{N}_i\cap\mathcal{C}_{p,k}\backslash\mathcal{M}_{p,k}, \label{eq:coup_con_coal}
\end{equation}
for all  $\forall i \in \mathcal{C}_{p,k}$ and $\ell \in\{k,\dots, k+h-1 \}$ (c.f. Remark \ref{rem:const_coal}), a distributed Lagrangian approach, where the coupling constraints \eqref{eq:coup_con_coal} are relaxed, can be implemented to solve Problem \eqref{eq:lo_mpc}. In this regard, a distributed dual-ascent algorithm, such as that presented in \cite{ananduta2018}, can be applied to solve Problem \eqref{eq:lo_mpc}, in which more than one microgrid is involved. Note that, different distributed algorithms, e.g., \cite{boyd2011,wang2015,kraning2014,baker2016}, can also be chosen instead. Nevertheless, since there are available distributed algorithms in the literature that can be applied, we assume that the optimal solution to Problem \eqref{eq:lo_mpc} can be computed in a distributed manner.

\section{Sub-optimality and communication trade-off}
\label{sec:opt}
In this section, we discuss the sub-optimality and communication trade-off of the proposed scheme. 
First, we show an estimation of the sub-optimality level achieved by the scheme. To that end, we state the collection of the optimization problems \eqref{eq:lo_mpc}, for all coalitions $\mathcal{C}_{p,k}, p=1,\dots,m$, as follows:
\begin{equation}
	\begin{aligned}
	&\min_{\{\{(\bm{u}_{i,\ell},\bm{v}_{i,\ell})\}_{i \in \mathcal{N}}\}_{\ell=k}^{k+h-1}} \sum_{i \in \mathcal{N}} \sum_{\ell =k}^{k+h-1}  J_{i,\ell}(\bm{u}_{i,\ell},\bm{v}_{i,\ell})\\
	&\text{s.t. \eqref{eq:co_loc_cons}, \eqref{eq:co_loc_cons2}, and \eqref{eq:coup_dec_const}, } \forall p \in \{1,\dots,m\},
	\end{aligned}
		\label{eq:lo_mpc_all}%
\end{equation}
for all $\ell \in\{k,\dots,k+h-1 \}$. Denote the optimal value of Problem \eqref{eq:lo_mpc_all} by $J^{\star}_k$. Note that $J^{\star}_k$ represents the cost function value of Problem \eqref{eq:MPC_net} computed by the proposed scheme. Furthermore, denote by $J_k^o$ the optimal value of Problem \eqref{eq:MPC_net} and define the sub-optimality measure as the difference between the cost function value computed using the proposed scheme and the optimal value of Problem \eqref{eq:MPC_net}, denoted by $\Delta J_k$, i.e., 
\begin{equation}
	\Delta J_k = J^{\star}_k-J^o_k. \label{eq:diff_J}
\end{equation}
\begin{prop}
	\label{prop:est_difJ}
	Let $J^{\star}_k$ and $J^o_k$ be the optimal values of Problems \eqref{eq:lo_mpc_all} and \eqref{eq:MPC_net} at time $k$, respectively. Furthermore, let $J^{b}_k$ denote the optimal value of the following optimization problem:
	\begin{equation}
	\begin{aligned}
	&\min_{\{\{(\bm{u}_{i,\ell},\bm{v}_{i,\ell})\}_{i \in \mathcal{N}}\}_{\ell=k}^{k+h-1}} \sum_{i \in \mathcal{N}} \sum_{\ell =k}^{k+h-1}  J_{i,\ell}(\bm{u}_{i,\ell},\bm{v}_{i,\ell})\\
	&\mathrm{s.t. } \text{ \eqref{eq:co_loc_cons} $\mathrm{and}$ \eqref{eq:co_loc_cons2} } \forall p \in \{1,\dots,m\}.
	\end{aligned}
	\label{eq:lb_mpc_all}%
	\end{equation}
	Then, the following estimate on the suboptimality measure $\Delta J_k$, defined in \eqref{eq:diff_J}, holds:
	\begin{equation}
		\Delta J_k \leq J^{\star}_k-J^b_k. \label{eq:est_diffJ}
	\end{equation}
\end{prop}
\begin{pf*}{Proof.}
	Note that Problem \eqref{eq:lb_mpc_all} can be obtained by relaxing Problem \eqref{eq:MPC_net}. In particular, the coupling constraint $v_{i,\ell}^j + v_{j,\ell}^i=0$, for each pair of nodes $i$ and $j$ that do not belong to the same coalition, is discarded in Problem \eqref{eq:lb_mpc_all}. Due to this relaxation, we can conclude that $J^b_k \leq J^o_k$. Moreover, based on Proposition \ref{prop:feas_of_orig}, the solution obtained by Problem \eqref{eq:lo_mpc}, for all coalitions $\mathcal{C}_{p,k}, \ p=1,\dots,m$, is also a feasible solution to Problem \eqref{eq:MPC_net}, implying that $J^o_k \leq J^{\star}_k$. Based on the preceding observations, the relation in \eqref{eq:est_diffJ} holds. \eod
\end{pf*}
\begin{rem}
	\label{re:tight_bound}
	Consider the case when $\mathcal{C}_{p,k}=\mathcal{N}$, for $p=1,\dots,m$. In this case, for any $i\in \mathcal{N}$, all neighbors of node $i$, i.e., $j\in\mathcal{N}_i$, belong to the same coalition as that of node $i$. Thus, in \eqref{eq:coup_dec_const},  $\mathcal{N}_i\backslash\mathcal{C}_{p,k}=\emptyset$. This fact implies that Problem \eqref{eq:lo_mpc_all} is equivalent to Problem \eqref{eq:MPC_net} and Problem \eqref{eq:lb_mpc_all}, implying   $\Delta J_k=0$ and $J^{\star}_k-J^b_k=0$. Additionally, Problem \eqref{eq:lb_mpc_all} can be decomposed into {$m$ sub-problems, each of which can be assigned to each coalition}.\eod
\end{rem}

Now, we discuss the communication cost of the proposed scheme. Algorithms \ref{alg:repartitioning} and \ref{alg:clust} do require information exchange among the controllers. The total size of data exchanged throughout the process in Algorithms \ref{alg:repartitioning} and \ref{alg:clust} is $\mathcal{O}(m)$ per iteration. 
Finally, we evaluate the size of data communicated when solving the coalition-based economic dispatch problem. Each coalition might need to use a distributed optimization method since there might be more than one microgrid in a coalition. As an example, we consider the dual-ascent algorithm \cite{ananduta2018} as the distributed optimization method. In this algorithm, the size of exchanged information is  $\mathcal{O}(m|\mathcal{N}|h)$ per iteration since each microgrid must exchange the coupled decision variables with each neighbor. 
In the best-case scenario, when all microgrids are self-sufficient, communication might not be necessary at all at one time instant. Furthermore, even if the repartitioning procedure is triggered, in the worst-case scenario, i.e., when the resulting coalition includes all microgrids, the extra amount of data must be exchanged to perform the repartitioning and coalition formation procedures is relatively much smaller than that of performing the distributed algorithm. In addition, for a coalition that has only one microgrid, its controller only needs to solve a local optimization problem once, which significantly reduces the computational burden as well.  
\color{black}

In comparison with existing methods that are based on distributed optimization algorithms, e.g., \cite{baker2016,wang2015,kraning2014,braun2016,guo2016}, the proposed scheme reduces the number of neighbors with which each agent must communicate since each agent only needs to communicate with a subset of neighbors that belong to the same coalition. This fact implies the reduction of communication flows. Moreover, the coalition-based problem \eqref{eq:lo_mpc} is relatively smaller than the network problem \eqref{eq:MPC_net}, thus intuitively a solution to \eqref{eq:lo_mpc} can be computed faster than a solution to \eqref{eq:MPC_net} using the same distributed iterative algorithm. 

Finally, we discuss the practicality of performing the proposed scheme. As in any distributed scheme, local controllers must cooperate to perform the proposed method and a communication network must also be available. Since the partition of the electrical network is dynamic, a dynamic network, containing necessary communication links, might be required. Another possibility is by having an all-to-all network although in the process not all links will be used. Furthermore, each local controller must also be able to communicate with the active components of the network, i.e., the storage and dispatchable generation units. The second important note that we would like to mention that although in this paper we consider an MPC-based framework, where the set points are computed at each time instant, the proposed method can also be implemented for day-ahead economic dispatch without requiring any modification. In this case, the prediction horizon is set to be one day and prior to the computation of the decisions, the self-sufficiency of each microgrid is evaluated. \color{black}

\section{Case study}
\label{sec:case_st}
\begin{figure}
	\centering
	\includegraphics[scale=0.9]{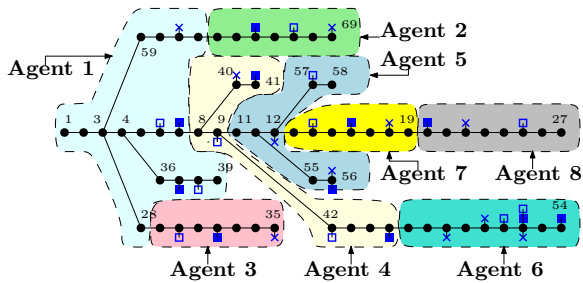}
	\caption{The topology of the PG\&E 69-bus distribution system and its 8-agent initial partition \cite{arefifar2012}. Squares indicate the distributed generation units, i.e., \textcolor{blue}{$\blacksquare$} and \textcolor{blue}{$\Box$} represent a renewable generation unit and a dispatchabale generator, respectively, whereas crosses, \textcolor{blue}{$\boldsymbol{\times}$}, indicate the storages.
	}
	\label{fig:cs_top}
\end{figure}
We consider the PG\&E 69-bus distribution network, as shown in Fig. \ref{fig:cs_top} where dispatchable, solar-based distributed generation, and storage units are added to the network. 
The simulation time is one day with the sampling time of 15 minutes, implying 96 time steps. Furthermore, the prediction horizon is set to be 8 time steps and the weight on the cost of the repartitioning problem is set as $\alpha=10^4$. 

\begin{figure}
	\centering
	\includegraphics[scale=0.6]{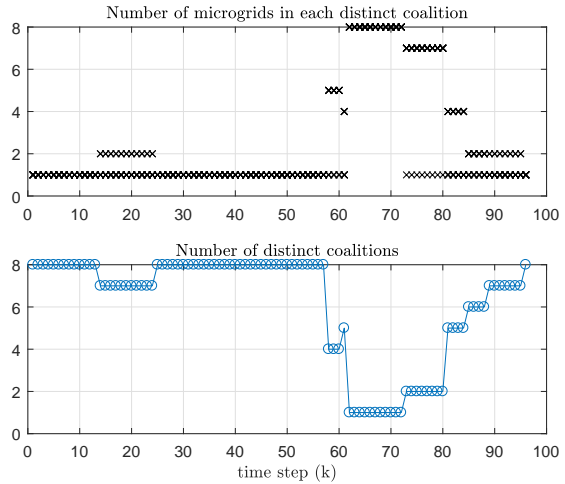}
	\caption{The evolution of coalitions formed.
	}
	\label{fig:co_res}
\end{figure}
\begin{figure}
	\centering
	\includegraphics[scale=0.6]{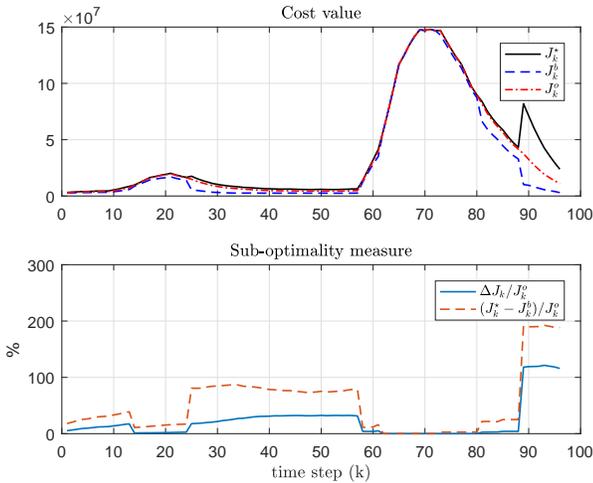}
	\caption{Top plot shows the cost values computed using the proposed scheme, $J^{\star}_k$, (solid line), by solving Problem \eqref{eq:MPC_net} centrally as the benchmark $J^{o}_k$, (dashed-dotted line), and the lower bound, $J^{b}_k$ (dashed line). Bottom plot shows the suboptimality of the proposed scheme and its upper bound.
	}
	\label{fig:deltaJ}
\end{figure}

The initial partition of the network is based on one of the partitioning results in \cite{arefifar2012}. How the microgrids form coalitions throughout the simulation can be seen in Figure \ref{fig:co_res}. 
We can observe that 
 during the peak hours $57\leq k \leq 80$, coalitions must be formed and, even at a certain period, all microgrids must join as one coalition, whereas during the off-peak hours, self-sufficient microgrids can be formed. 
Figure \ref{fig:deltaJ} shows the stage costs for all time steps and the sub-optimality of the proposed scheme. 
\color{black}


\section{Conclusion and future work}
\label{sec:concl}
We develop a non-centralized MPC-based economic dispatch scheme for systems of interconnected microgrids. The approach consists of an event-triggered repartitioning method with the aim of maintaining self-sufficiency of each microgrid and decomposing the centralized economic dispatch problem into coalition-based sub-problems in order to compute a feasible but possibly sub-optimal decisions. The main advantage of the approach is a low communication burden, which is essential for online applications. 
\vspace{-10pt}
\begin{ack}                               
This work has received funding from the European Union's Horizon 2020 research and innovation programme under the Marie Sk\l{}odowska-Curie grant agreement No 675318 (INCITE).  
\end{ack}
\vspace{-10pt}
{

\bibliographystyle{dcu}        
\bibliography{ref_ja}           
}                                

\end{document}